\documentstyle[11pt]{article}
\newcommand{\bbz}{\mbox{\boldmath $Z$}}
\newcommand{\bbn}{\mbox{\boldmath $N$}}
\newcommand{\bbc}{\mbox{\boldmath $C$}}
\newcommand{\bbr}{\mbox{\boldmath $R$}}
\newcommand{\bbf}{\mbox{\boldmath $F$}}
\newcommand{\qed}{\hfill$\Box$}
\newcommand{\major}{\sl maj}
\newcommand{\var}{\hbox{Var}}
\newcommand{\descent}{\sl des}
\newcommand{\sh}{\sl shape}
\newcommand{\pr}{\hbox{Pr}}
\newcommand{\dt}{\hbox{det}}
\newenvironment{note}[1]{\par\addvspace{\medskipamount}\noindent
                         {\bf {#1}}\sl
                       }{\par\addvspace{\medskipamount}\rm}

\title{Descent Functions and Random Young Tableaux}
\author{Ron M. Adin%
\thanks{Department of Mathematics and Computer Science, Bar-Ilan University,
Ramat-Gan 52900, Israel. Email: {\tt radin@math.biu.ac.il}}
$^\ddag$
\and Yuval Roichman%
\thanks{Department of Mathematics and Computer Science, Bar-Ilan University,
Ramat-Gan 52900, Israel. Email: {\tt yuvalr@math.biu.ac.il}}
\thanks{Both authors supported in part by the Israel Science Foundation
and by internal research grants from 
Bar-Ilan University.}}
\date{submitted October 5, 1999\ ;\ revised October 12, 2000}

\begin{document}

\maketitle
\begin {abstract}
The expectation of the descent number of a random Young tableau of a fixed shape is given,
and concentration around the mean is shown.
This result is generalized to the major index and to other descent functions.
The proof combines probabilistic arguments 
 together with combinatorial character theory.
Connections with Hecke algebras are mentioned.

\end{abstract}

\section{\bf Introduction}

\subsection{Background}

In the late sixties Erd\H os and Tur\'an 
have published a classical series of papers
on random permutations.
Since then there has been a resurgence of interest in probabilistic aspects of combinatorial parameters of permutations and related objects. 

This paper deals principally with two classical
combinatorial parameters: descent number and major index.
These parameters were originally studied in the context of permutations.
The study of the descent number of a permutation started with Euler;
the major index has been introduced by MacMahon [M].
Foata-Sch\"utzenberger [F, FS], Garsia-Gessel [GG] and others
carried out an extensive research of these parameters.
The definitions of descent number and major index for permutations lead to definitions
of the same concepts for Young tableaux.
These parameters on permutations -- as well as on tableaux -- play significant roles
in algebraic combinatorics : the Solomon descent algebra [Re, Ch. 9],
Schur functions [St, Ch. 7], and combinatorial character formulas [Ro2].
These concepts were
 also applied to sorting [K, Section 5.1] and card shuffling [DMP].

In this paper we study the distribution of these parameters for random Young tableaux of a given shape. 
Proofs of the main results are obtained by 
a combination of probabilistic arguments and combinatorial character theory.

\subsection{Main Results}

Let $\lambda$ be a partition of $n$ 
(For definitions of basic concepts see Section 2).
We shall be concerned with {\it random} (standard Young) tableaux,
assumed to be chosen uniformly with prescribed shape $\lambda$.
A {\it descent} in a standard Young tableau $T$ is an entry
$i$ such that $i+1$ is strictly south (and weakly west) of $i$.
Denote the set of all descents in $T$ by $D(T)$.

For any 
function $f:\bbn\rightarrow \bbr$ we introduce
the corresponding {\it descent function} $d_f$ on standard Young tableaux :
$$
d_f(T):=\sum\limits_{i\in D(T)} f(i).
$$
This concept generalizes two classical combinatorial parameters, the
{\it descent number} and the {\it major index} (for tableaux) :
$$
\descent(T):=\sum\limits_{i\in D(T)} 1\hbox{ };\qquad
\major(T):=\sum\limits_{i\in D(T)} i.
$$

\medskip

In this paper we prove

\begin{note}
\noindent{\bf Theorem 1.} Let $\lambda$ be a fixed partition of $n$,
 and let $E_\lambda[d_f]$ be the expected value of 
a descent function $d_f$ on random standard Young 
tableaux of shape $\lambda$.
Then
$$
E_\lambda [d_f]= c(\lambda)\cdot\sum\limits_{i=1}^{n-1}f(i).
$$
where $c(\lambda):= [{n\choose 2}-\sum\limits_i{\lambda_i\choose 2}+
\sum\limits_j{\lambda_j'\choose 2}]/ n(n-1)$.
(Here $\lambda_i$ is the length of the $i$-th row in the Young diagram of shape $\lambda$, and $\lambda_i'$ is the length of the $i$-th column.)
\end{note}

\smallskip

\noindent See Theorem 4.1 below.

\smallskip

Under mild conditions the descent function is concentrated around its mean.
A function $f:\bbn\rightarrow \bbr$ has {\it strictly polynomial growth}
if there exist constants $0<c_1<c_2$ and $\alpha>0$, such that
 $c_1\le {f(n)\over n^\alpha}\le c_2$ for $n$ large enough.

\smallskip

\begin{note}
\noindent{\bf Theorem 2.}
Let $0<\delta<1$ and $\varepsilon>0$ be fixed constants, and
let $\lambda$ be a partition of $n$ with $\lambda_1\le \delta n$.
  Then for any function $f$ with strictly polynomial growth,
and for a random standard Young tableau $T$
of shape $\lambda$,
$$
d_f(T)=(1+O(n^{-{1\over 2}+\varepsilon}))E_\lambda[d_f] 
$$
almost surely (i.e., with probability tending to 1 as $n$ tends to infinity),
uniformly on $f$ and $\lambda$ as above.
\end{note}

\smallskip

\noindent See Theorem 5.1 below.

\smallskip

 For other work, following the current paper,  see [H].

\medskip

The rest of the paper is organized as follows.
Definitions, notations and necessary preliminaries are given
in Section 2.
In Section 3 the expectation and variance of the major index are
evaluated. Here Stanley's hook formula plays a crucial role.
Results obtained in Section 3 are extended to general descent functions in Sections 4 and 5 by combining probabilistic arguments with combinatorial
character formulas.  We end the paper with remarks on connections
of the statistics of descent functions with the spectra of so called `good'
elements in Hecke algebras.

\medskip

\section{Preliminaries}

\subsection{Young Tableaux}

Let $n$ be a positive integer. A {\it partition} of $n$
is a vector of positive integers
$\lambda=(\lambda_1,\lambda_2,\ldots,\lambda_k)$,
where 
$\lambda_1\ge\lambda_2\ge\ldots\ge\lambda_k$
and  $\lambda_1+\ldots+\lambda_k=n$.   
For a partition  $\lambda=(\lambda_1,\ldots,\lambda_k)$  define the {\it conjugate partition}
 $\lambda'=(\lambda'_1,\dots,\lambda'_t)$ by letting $\lambda'_i$ be the number of parts of $\lambda$ that are $\ge i$. 

\smallskip

The {\it dominance (partial) order} on partitions is defined as follows :        
For any two partitions of $n$, $\mu$ and $\lambda$ ,
$\mu$ {\it dominates} $\lambda$ 
if and only if for any $i$ $\ $
$\sum_{j=0}^i \lambda_j\le \sum_{j=0}^i \mu_j.$

\smallskip

 For example, $\lambda=(4,4,2,1)$ is a partition of 11. Then  $\lambda'=(4,3,2,2)$, and $\lambda$ dominates $\lambda'$.

\medskip

The set 
 $\{(i,j)\ |\ i,j\in \bbz$ , $ 0<i\le k$ , $0<j\le\lambda_i\}$ is called the {\it Young diagram}
of    {\it shape} $\lambda$.
$(i,j)$ is the {\it cell} in row $i$ and column $j$.
The   diagram of the conjugate
shape $\lambda'$ may be obtained from the diagram of shape $\lambda$ 
by interchanging rows and columns.

\medskip

  A {\it Young tableau} of shape $\lambda$
is obtained by inserting the integers $1,2,\ldots,n$
as {\it entries} in the
cells of the Young diagram of shape $\lambda$, allowing no repetitions.
A {\it standard Young tableau} of shape $\lambda$ is a 
Young tableau whose entries increase along rows and columns.

\smallskip

We shall draw Young tableaux as in the following example.

\noindent{\bf Example 1.}

\[
\begin{array}{ccccc}
1         & 3         & 4                 & 6     & 9     \\
2         & 7         & 8                 & 12        \\
5         & 11          \\
10             \\
\end{array}
\]

\medskip

The {\it hook length} of a cell $(i,j)$ in the diagram of shape $\lambda$ is defined by                                            
$$
h_{i,j} := \lambda_i+\lambda_j'-i-j+1.
$$

Denote by  $f^\lambda$ the number of standard  Young tableaux
of shape $\lambda$. 
A famous combinatorial formula describes this number 
in term of hook lengths.
      
\smallskip

\begin{note}
\noindent{\bf The Frame-Robinson-Thrall Hook Formula}. [Sa, Theorem 3.1.2]
$$
f^\lambda=  {n!\over \prod\limits_{(i,j)\in\lambda}h_{ij}}.
$$
\end{note}

\bigskip

\subsection{Descents}

A {\it descent} in a standard Young tableau $T$ is an entry
$i$ such that $i+1$ is strictly south (and weakly west) of $i$.
Denote the set of all descents in $T$ by $D(T)$.
The {\it descent number} and the {\it major index} (for tableaux) 
are defined as follows :
$$
\descent(T):=\sum\limits_{i\in D(T)} 1\hbox{ };\qquad
\major(T):=\sum\limits_{i\in D(T)} i.
$$

\noindent{\bf Example 2.}
Let $T$ be the standard Young tableau drawn in Example 1.
Then $D(T)=\{1,4,6,9\}$, $\descent(T)=4$, and $\major(T)=1+4+6+9=20$.

\medskip

The following theorem describes the generating function for the
major index of standard Young tableaux.

\medskip

\begin{note} 
\noindent{\bf The Stanley Hook Formula} [St, Corollary 21.5]
$$
\sum_{{\sl shape}(T)=\lambda} q^{\major(T)}=q^{\sum_i (i-1)\lambda_i} \cdot {\prod\limits_{k=1}^n [k]_q\over \prod\limits_{(i,j)\in \lambda} [h_{ij}]_q},
$$
where the sum is taken over all standard Young tableaux of shape $\lambda$,
$h_{ij}$ are the hook lengths in the diagram of $\lambda$,
and for any positive integer $m$
$$
[m]_q:=1+q+q^2+\dots+q^{m-1}.
$$
\end{note}

\medskip

For $q=1$ this formula reduces to the Frame-Robinson-Thrall 
hook formula for the number of standard Young tableaux of  a given shape. 
No such formula is known for the descent number of tableaux.

\subsection{Characters}

A (complex) representation of a group $G$ is  
a homeomorphism $\rho:G\rightarrow GL_n(\bbc)$.
The character $\chi^\rho:G\rightarrow \bbc$ is the trace
of $\rho(g)$, $g\in G$. 
By definition, the character is a class
function on the group (i.e., invariant under conjugation).
An {\it irreducible} representation is a representation
which has no nontrivial subspace invariant under all $\rho(g)$, $g\in G$.

The conjugacy classes of the symmetric group $S_n$
are described by their cycle type; thus, by the partitions of $n$ .
The irreducible representations of $S_n$ are also indexed
by these partitions. See e.g. [Sa].

Let $\lambda$ and $\mu$ be partitions of $n$.
Denote by $\chi^\lambda_\mu$ the value at a conjugacy class
of cycle type $\mu$
of the character of 
the irreducible representation indexed by $\lambda$.

\medskip

The following combinatorial formula represents the irreducible
characters of $S_n$ in terms of descents of standard Young tableaux.
This formula is a special case of [Ro1, Theorem 4].

\smallskip

\begin{note}
\noindent{\bf Theorem 2.1.} 
$$
\chi^\lambda_\mu=\sum\limits_{\sh(T)=\lambda} \hbox{weight}_\mu(T),
$$
where the sum is taken over all standard tableaux of shape $\lambda$,
and the weight 
$\hbox{weight}_\mu(T)\in \{\pm 1, 0\}$ is defined as 
follows :
$$
\hbox{weight}_\mu(T) := 
\prod_{1\le i\le k\atop i\not\in B(\mu)} f_\mu(i,T),
$$
where 
$B(\mu)=\{\mu_1, \mu_1+\mu_2,\dots,\mu_1+\cdots+\mu_t\}$,
and
$$
f_\mu(i,T) := 
\cases
{-1 &  $i\in D(T)$;\cr
0 &  $i\not\in D(T), i+1\in D(T)$  and   $i+1\not\in B(\mu)$;\cr
1 &  otherwise. \cr}
$$
\end{note}

\medskip

\section{\bf Major Index}

 In this section we apply Stanley's hook formula 
(see Subsection 1.2) to evaluate the 
expectation and variance of the major index of random Young tableaux
of a given shape.

\medskip

\begin{note}
\noindent{\bf Proposition 3.1}
Let $\lambda$ be a fixed partition of $n$,
 and let $E_\lambda[\major]$ be the expected value of the major index on random standard Young 
tableaux of shape $\lambda$.
Then
$$
E_\lambda [\major]= 
{1\over 2}\left[{n\choose 2}-
\sum_i {\lambda_i\choose 2}
+\sum_j {\lambda'_j\choose 2}\right].
$$
\end{note}

\smallskip

\noindent{\bf Proof.} Let 
$$
f(q):=\sum_{\sh(T)=\lambda} q^{\major(T)}.
$$
Then
$$
q \cdot f'(q)= \sum_{\sh(T)=\lambda} \major (T) \cdot q^{\major(T)}.
$$
Hence
$$
(\log f(q))'|_{q=1}=\left.{f'(q)\over f(q)}\right|_{q=1}=
{ \sum_{\sh(T)=\lambda} \major (T) \over  \sum_{\sh(T)=\lambda} 1}=
E_\lambda[\major]. \leqno (3.1)
$$

\smallskip

In order to evaluate the expected value of the
major index we need the following elementary limit :
$$
\lim_{q\rightarrow 1} {[m]_q' \over [m]_q}={{m\choose 2}\over m}={m-1\over 2}. \leqno (3.2)
$$

\smallskip

Substituting Stanley's hook formula in (3.1), and using (3.2)
we obtain
$$
E_\lambda[\major]=(\log f(q))'|_{q=1}=
\sum_i (i-1)\lambda_i + \sum_{k=1}^n {k-1\over 2} - \sum_{(i,j)\in \lambda}
{h_{ij}-1\over 2}. \leqno (3.3)
$$
Note that
$$
\sum_i (i-1)\lambda_i=\sum_j {\lambda'_j\choose 2} \leqno (3.4)
$$
and
$$
\sum_{(i,j)\in \lambda}{h_{ij}-1\over 2}={1\over 2} \left[\sum_i 
{\lambda_i\choose 2}+\sum_j{\lambda_j'\choose 2}\right]. \leqno (3.5)
$$
Substituting (3.4) and (3.5) in (3.3) completes the proof.

\qed

\medskip

\begin{note}
\noindent{\bf Proposition 3.2}
Let $\lambda$ be a fixed partition of $n$,
 and let $\var_\lambda[\major]$ be the variance of the major index on random standard Young 
tableaux of shape $\lambda$.
Then
$$
\var_\lambda [\major]= {1\over 12}\left[\sum_{k=1}^n k^2 - \sum_{(i,j)\in \lambda} h_{ij}^2\right].
$$
\end{note}

\smallskip

\noindent{\bf Proof.} Let $f$ denote, as in the previous proof,
the generating function of the major index of Young tableaux of a given shape. Then
$$
\sum_{\sh(T)=\lambda} (\major (T))^2\cdot q^{\major(T)}=
q^2 \cdot f''(q) + q \cdot f'(q).
$$
Hence
$$
\var_\lambda[\major]= E_\lambda[\major^2]-
(E_\lambda[\major])^2 = \leqno (3.6)
$$
$$
=\left.\left[{q^2 f'' + q f'\over f} - \left(q{f'\over f}\right)^2\right]\right|_{q=1} = 
\left.\left[q^2 \left({f'\over f}\right)' + q{f'\over f}\right]\right|_{q=1}.
$$
Now
$$
{f'\over f}= (\log f)'=
\sum_i (i-1)\lambda_i \cdot q^{-1}+
\sum_{k=1}^n {[k]_q'\over [k]_q}-
\sum_{(i,j)\in\lambda} {[h_{ij}]'_q\over [h_{ij}]_q}
$$
and
$$
\lim_{q\rightarrow 1}\left({[m]_q'\over [m]_q}\right)'=
\lim_{q\rightarrow 1}\left[ {[m]_q''\over [m]_q}-
\left({[m]_q'\over [m]_q}\right)^2\right] =
$$
$$
=\lim_{q\rightarrow 1}\left[ {\sum_{k=2}^m(k-1)(k-2)q^{k-3}\over
\sum_{k=1}^{m}q^{k-1}}-
\left({\sum_{k=1}^m(k-1)q^{k-2}\over \sum_{k=1}^m q^{k-1}}\right)^2\right]=
$$
$$
= {\sum_{k=2}^m(k-1)(k-2)\over m}-\left({\sum_{k=1}^m(k-1)\over m}\right)^2=
{2{m\choose 3}\over m}-\left({m-1\over 2}\right)^2={(m-1)(m-5)\over 12}.
$$
Hence
$$
\lim_{q\rightarrow 1}\left({f'\over f}\right)'=
-\sum_i(i-1)\lambda_i+\sum_{k=1}^n {(k-1)(k-5)\over 12}-
\sum_{(i,j)\in \lambda} {(h_{ij}-1)(h_{ij}-5)\over 12}.\leqno(3.7)
$$
Substituting (3.7) and (3.3) into 
the right hand side of (3.6) we obtain the desired result.

\qed

\medskip

\begin{note}
\noindent{\bf Corollary 3.3}
Let $0<\delta<1$ and $\varepsilon>0$ be fixed constants, and
let $\lambda$ be a partition of $n$ with $\lambda_1\le \delta n$.
Then
$$
\major(T)=(1+O(n^{-{1\over 2}+\varepsilon}))E_\lambda[\major] 
$$
almost surely (i.e., with probability tending to 1 as $n$ tends to infinity)
uniformly on $\lambda$ as above.
\end{note}

\smallskip

\noindent{\bf Proof.} By Chebyshev's inequality [Fe, (6.2)] 
$$
\hbox{Pr}\left(|\major(T)-E_\lambda[\major]|\ge t E_\lambda[\major]\right)
\le {1\over t^2}\cdot {\var_\lambda[\major]\over E_\lambda[\major]^2}.
$$

Now, by Proposition 3.2
$$
0\le \var_\lambda[\major]\le {1\over 12} \sum_{k=1}^n k^2 = O(n^3).
$$

In order to bound $E_\lambda[\major]$ from below, note that
$$
\sum_i{\lambda_i\choose 2}-\sum_j{\lambda_j'\choose 2}=\sum_{(i,j)\in\lambda} (j-i). \leqno (3.8)
$$ 
Thus this expression is a monotone increasing function of $\lambda$ with respect to the dominance order of partitions (see also [Su]).
So, under the restriction $\lambda_1\le \delta n$ (we may assume that
$\delta> {1\over 2}$) this expression is maximized
when $\lambda_1=\delta n, \lambda_2=(1-\delta)n$. Hence, by Proposition 3.1 :
$$
E_\lambda[\major]=
{1\over 2}\left[{n\choose 2}
-\sum_i{\lambda_i\choose 2}+\sum_j{\lambda_j'\choose 2}\right]
={1\over 2}\left[{n\choose 2}-\sum_{(i,j)\in \lambda}(j-i) \right]\ge
$$
$$
\ge {1\over 2}\left[{n\choose 2}-{\delta n\choose 2}-
{(1-\delta)n\choose 2}\right]=\Omega(n^2).
$$
Therefore,
$$
{1\over t^2}\cdot {\var_\lambda[\major]\over E_\lambda[\major]^2}\le
O\left({1\over t^2 n}\right)=O(n^{-2\varepsilon}),
$$
provided that ${1\over t}=O(n^{1/2-\varepsilon})$.
 We conclude that 
for such $\lambda$ and $t$
$$
(1-t)E_\lambda[\major] < \major(T)< (1+t)E_\lambda[\major] 
$$
with probability tending to 1 as $n$ tends to infinity.

\qed

\bigskip

\section{Expectation of Descent Functions}

In this section we generalize Proposition 3.1 to an arbitrary descent function.
This is done by combining probabilistic arguments 
 together with combinatorial character formulas.

Recall the definition of descent functions from Section 1.

\medskip

\begin{note}
\noindent{\bf Theorem 4.1} Let $\lambda$ be a fixed partition of $n$,
let $f:\bbn\rightarrow \bbr$ be an arbitrary function,
 and let $E_\lambda[d_f]$ be the expected value of $d_f$ on random standard Young 
tableaux of shape $\lambda$.
Then
$$
E_\lambda [d_f]= c(\lambda)\cdot\sum\limits_{i=1}^{n-1}f(i).
$$
where $c(\lambda):= [{n\choose 2}-\sum\limits_i{\lambda_i\choose 2}+
\sum\limits_j{\lambda_j'\choose 2}]/ n(n-1)$.
\end{note}

\medskip

Substituting $f(i)=1$ and $f(i)=i$ in Theorem 4.1 gives the expectations
for the descent number and for the major index of a random tableau, respectively. In particular, Proposition 3.1 is a special case of Theorem 4.1.

\bigskip

To prove Theorems 4.1 and 5.1 we shall use 
a variant of the character formula
given in Theorem 2.1.

\medskip

\noindent For a fixed $1\le i < n$ define the $i$-weight 
$$
\hbox{weight}^i_{(21\dots 1)}(T):=\cases
{-1, & if $i\in D(T)$ ;  \cr
1, & if $i\not\in D(T)$ .   \cr}
$$
For a fixed $1\le i < n-1$ define another $i$-weight
$$
\hbox{weight}^i_{(31\dots 1)}(T):=\cases
{-1, & if $i\in D(T)$ and $i+1\not\in D(T)$  ; \cr
1, & if $i,i+1\in D(T)$ or $i,i+1\not\in D(T)$ ;  \cr
0, & if $i\not\in D(T)$ and $i+1\in D(T)$ .  \cr}
$$
For a fixed pair $1\le i<j<n$ with $j-i>1$ define the $ij$-weight
$$
\hbox{weight}^{ij}_{(221\dots 1)}(T):=\cases
{1, & if $i,j\in D(T)$ or $i,j\not\in D(T)$  ; \cr
-1, & otherwise  . \cr}
$$ 
The $i$- (or $ij$-) weight of a standard tableau $T$ depends on $i$
(or $i,j$). However, the sum of $i$-weights (or $ij$-weights) over
all standard tableaux of shape $\lambda$ is independent of $i$ and $j$
and gives the corresponding character:

\smallskip

\begin{note}

\noindent{\bf Lemma 4.2} For any partition $\lambda$ of $n$,
$$
\chi^\lambda_{(21\dots 1)}=\sum\limits_{\sh(T)=\lambda} \hbox{weight}^i_{(21\dots 1)} (T) \qquad\qquad (1\le i < n) \leqno \hbox{(i)}
$$
$$
\chi^\lambda_{(31\dots 1)}=\sum\limits_{\sh(T)=\lambda} \hbox{weight}^i_{(31\dots 1)} (T) \qquad\qquad (1\le i < n-1) \leqno \hbox{(ii)}
$$
$$
\chi^\lambda_{(221\dots 1)}=\sum\limits_{\sh(T)=\lambda} \hbox{weight}^{ij}_{(221\dots 1)} (T) \qquad\qquad (1\le i< j-1 < n-1) 
\leqno \hbox{(iii)}
$$

\end{note}

\smallskip

\noindent For proofs and more details see [Ro1, Section 7].

\medskip

Recall also that
$$
\chi^\lambda_{(1\dots 1)} =\sum\limits_{\sh(T)=\lambda} 1 =f^\lambda.
$$

\medskip

\noindent{\bf Proof of Theorem 4.1.} 
Let $T$ be a random standard Young tableau of shape $\lambda$.
For $1\le i< n$,
let $X_i$ be the random variable defined by
$$
X_i:=\cases
{1, & if $i\in D(T)$ ; \cr
0, & otherwise  . \cr}
$$
Then
$$
E_\lambda [d_f]=
E_\lambda\left[\sum\limits_{i=1}^{n-1}f(i)X_i\right]=
\sum\limits_{i=1}^{n-1}f(i)E_\lambda[X_i]. \leqno (4.1)
$$
Now, by definition
$$
\hbox{weight}^i_{(21\dots 1)}=1-2X_i
$$
and therefore, by Lemma 4.2(i)
$$
1-2E_\lambda[X_i]=E_\lambda[\hbox{weight}^i_{(21\dots 1)}]=
\chi^\lambda_{(21\dots 1)}/ \chi^\lambda_{(1\dots 1)} . \leqno (4.2)
$$
In particular, note that $E_\lambda[X_i]$ is independent of $i$
(See [St, Prop. 7.19.9]).

A classical formula of Frobenius shows that 
$$
\chi^\lambda_{(21\dots 1)}/ \chi^\lambda_{(1\dots 1)}={1\over {n\choose 2}}
\left[\sum\limits_i {\lambda_i\choose 2}-\sum\limits_j{\lambda_j'\choose 2}\right]. \leqno (4.3)
$$ 
See [I].

Combining (4.1), (4.2) and (4.3) completes the proof.

\qed

\smallskip

\section{Concentration}

In this section Corollary 3.3 is generalized to descent functions, satisfying certain mild conditions.

Recall that a function $f:\bbn\rightarrow \bbr$ has {\it strictly polynomial growth}
if there exist constants $0<c_1<c_2$ and $\alpha>0$, such that
 $c_1\le {f(n)\over n^\alpha}\le c_2$ for $n$ large enough.

\smallskip

\begin{note}
\noindent{\bf Theorem 5.1}
Let $0<\delta<1$ and $\varepsilon>0$ be fixed constants, and
let $\lambda$ be a partition of $n$ with $\lambda_1\le \delta n$.
  Then for any function $f$ with strictly polynomial growth,
and for a random standard Young tableau $T$
of shape $\lambda$,
$$
d_f(T)=(1+O(n^{-{1\over 2}+\varepsilon}))E_\lambda[d_f] 
$$
almost surely (i.e., with probability tending to 1 as $n$ tends to infinity),
uniformly on $f$ and $\lambda$ as above.
\end{note}

\smallskip

\noindent This holds, in particular, if $d_f$ is either descent number or major index.

\smallskip

Theorem 5.1 is proved by estimating the asymptotic behavior of the variance.
This is done by 
expressing the variance in terms of $S_n$-characters evaluated 
at ``small" conjugacy classes, and 
showing that ${\hbox{Var}_\lambda[d_f]/ E_\lambda[d_f]^2}$ is independent 
of the function
$f$, up to a multiplicative constant.

\medskip

Denote
$$
r^\lambda_2:={\chi^\lambda_{21\dots 1}/ \chi^\lambda_{1\dots 1}} ,
$$
$$
 r^\lambda_3:={\chi^\lambda_{31\dots 1}/ \chi^\lambda_{1\dots 1}} ,
$$ 
$$
r^\lambda_{22}:={\chi^\lambda_{221\dots 1}/ \chi^\lambda_{1\dots 1}} .
$$ 
These are the values of the normalized irreducible character
corresponding to $\lambda$ at conjugacy classes of types $(21\dots 1), (31\dots 1)$
and $(221\dots 1)$ respectively.

\smallskip

The following lemma plays a crucial role in the proof of Theorem 5.1.

\smallskip

\begin{note}

\noindent{\bf Lemma 5.2.}
$$
\pr[X_i=1]={1\over 2}(1-r^\lambda_2) \qquad (1\le i< n). \leqno \hbox{(i)}
$$
$$
\pr[X_i=1 \hbox{ and } X_{i+1}=1]=
{1\over 2}(1-r^\lambda_2)-{1\over 3}(1-r^\lambda_3)
\qquad (1\le i < n-1). \leqno \hbox{(ii)}
$$
$$
\pr[X_i=1 \hbox{ and } X_j=1]={1\over 2}(1-r^\lambda_2)-
{1\over 4}(1-r^\lambda_{22})
\qquad (1\le i< j-1<n-1), \leqno \hbox{(iii)} 
$$
where the probability $\pr[\cdot]$ is taken in the probability
space of all standard Young tableaux of a given shape, defined in Section 1.

\end{note}

\smallskip

\noindent{\bf Proof.}

\noindent {\bf (i)} Since $X_i$ is a 0-1 variable,
$$
E_\lambda[X_i]=\pr[X_i=1].
$$
Using equation (4.2) and the definition of $r^\lambda_2$,
the desired result follows.

\bigskip

\noindent {\bf (ii)} By Lemma 4.2(ii), for $1\le i< n-1$ :
$$
r^\lambda_3=
\pr[X_i=1\wedge X_{i+1}=1]+\pr[X_i=0\wedge X_{i+1}=0]
-\pr[X_i=1\wedge X_{i+1}=0].
$$
Therefore, 
$$
r^\lambda_3=(1-\pr[X_i=1\wedge X_{i+1}=0]-
\pr[X_i=0\wedge X_{i+1}=1])-\pr[X_i=1\wedge X_{i+1}=0].
\leqno(5.1)
$$
Now
$$
\pr[X_i=1\wedge X_{i+1}=0]=\pr[X_i=1]-\pr[X_i=1\wedge X_{i+1}=1],
$$
$$
\pr[X_i=0\wedge X_{i+1}=1]=\pr[X_{i+1}=1]-\pr[X_i=1\wedge X_{i+1}=1],
$$
and therefore, using (i) above
$$
\pr[X_i=1\wedge X_{i+1}=0]=
\pr[X_i=0\wedge X_{i+1}=1]=
$$
$$
={1\over 2}(1-r^\lambda_2)-\pr[X_i=1\wedge X_{i+1}=1].
$$
Thus, from (5.1):
$$
1-r^\lambda_3=2\pr[X_i=1\wedge X_{i+1}=0]+\pr[X_i=0\wedge X_{i+1}=1]=
$$
$$
={3\over 2}(1-r^\lambda_2)-3\pr[X_i=1\wedge X_{i+1}=1],
$$
and (ii) follows.

\bigskip

\noindent {\bf (iii)} By Lemma 4.2(iii), for any $1\le i< j-1< n-1$ 
$$
r^\lambda_{22}=\pr[X_i=1\wedge X_j=1]+\pr[X_i=0\wedge X_j=0]-
$$
$$
-\pr[X_i=0\wedge X_j=1]-\pr[X_i=1\wedge X_j=0].
$$
Continue as in the proof of (ii) above.

\qed

\bigskip

\noindent Denote
$$
P_{2,\lambda}:=\pr[X_i=1] ,
$$
$$
P_{3,\lambda}:=\pr[X_i=1 \hbox{ and } X_{i+1}=1] ,
$$
$$
P_{22,\lambda}:=\pr[X_i=1 \hbox{ and } X_j=1] \qquad (j-i>1).
$$
By Lemma 5.2, these probabilities are well defined (independent of $i$
and $j$).

\bigskip

\noindent{\bf Proof of Theorem 5.1.} 
By Chebyshev's inequality
$$
\hbox{Pr}\left(|d_f-E_\lambda[d_f]|\ge t E_\lambda[d_f]\right)
\le {1\over t^2}\cdot {\var_\lambda[d_f]\over E_\lambda[d_f]^2}.
$$
In order to prove Theorem 5.1 it suffices to give an effective upper bound
on ${\var_\lambda[d_f]/ E_\lambda[d_f]^2}$.

\smallskip

For a random tableau $T$ of shape $\lambda$, 
$$
d_f(T)=\sum_{i=1}^{n-1} f(i)X_i.
$$
Thus,
$$
\var_\lambda[d_f]=E_\lambda[d_f^2]-E_\lambda[d_f]^2=
E_\lambda\left[(\sum_{i=1}^{n-1}X_i f(i))^2\right]-
\left(E_\lambda\left[\sum_{i=1}^{n-1}X_i f(i)\right]\right)^2= \leqno (5.2)
$$
$$
=\sum_{i=1}^{n-1} E_\lambda[X_i^2] f(i)^2+
2\sum_{i=1}^{n-2} E_\lambda [X_i X_{i+1}] f(i) f(i+1)+
$$
$$
+ 2\sum_{j-i>1} E_\lambda [X_i X_j] f(i) f(j)
- (\sum_{i=1}^{n-1} E_\lambda [X_i] f(i) )^2 .
$$

\medskip

Since $X_i$ is a $0-1$ variable
$$
E_\lambda[X_i^2]=E_\lambda[X_i]=P_{2,\lambda} \qquad (1\le i < n),
\leqno (5.3)
$$
$$
E_\lambda[X_iX_{i+1}]=P_{3,\lambda} \qquad (1\le i < n-1),
\leqno (5.4)
$$
and
$$
E_\lambda[X_iX_j]=P_{22,\lambda} \qquad (j-i>1) . \leqno (5.5)
$$
Substituting (5.3)-(5.5) into the right hand side of (5.2) we obtain
$$
\var_\lambda[d_f]=
P_{2,\lambda} \sum_{i=1}^{n-1} f(i)^2+
2P_{3,\lambda}\sum_{i=1}^{n-2} f(i) f(i+1)+
$$
$$
+ 2P_{22,\lambda}\sum_{j-i>1} f(i) f(j) -
P_{2,\lambda}^2 (\sum_{i=1}^{n-1} f(i) )^2.
$$
Denote 
$$
\Sigma_1:=\sum_{i=1}^{n-1}f(i),
$$
$$
\Sigma_2:=\sum_{i=1}^{n-1}f(i)^2,
$$ 
$$
\Sigma_3:=2\sum_{i=1}^{n-2}f(i)f(i+1).
$$
Obviously
$$
2\sum_{j-i>1} f(i) f(j)=\Sigma_1^2-\Sigma_2-\Sigma_3,
$$
and therefore
$$
\var_\lambda[d_f]=
P_{2,\lambda}\Sigma_2+
P_{3,\lambda}\Sigma_3+
P_{22,\lambda}(\Sigma_1^2-\Sigma_2-\Sigma_3)
-P_{2,\lambda}^2\Sigma_1^2=
$$
$$
=(P_{2,\lambda}-P_{22,\lambda})\Sigma_2+
(P_{3,\lambda}-P_{22,\lambda})\Sigma_3+
(P_{22,\lambda}-P_{2,\lambda}^2)\Sigma_1^2.
$$
Also
$$
E_\lambda[d_f]=P_{2,\lambda}\Sigma_1,
$$
and consequently, using Lemma 5.2(i)-(iii),
$$
{\var_\lambda[d_f]\over E_\lambda[d_f]^2}
={1-r^\lambda_{22}\over (1-r^\lambda_2)^2}\cdot
{\Sigma_2\over \Sigma_1^2}+
{4r^\lambda_3-3r^\lambda_{22}-1\over 3(1-r^\lambda_2)^2}\cdot
{\Sigma_3\over \Sigma_1^2}+
{r^\lambda_{22}-(r^\lambda_2)^2\over (1-r^\lambda_2)^2}. \leqno (5.6)
$$

\medskip

$f$ has a strictly polynomial growth. It follows that
$$
{\Sigma_2\over \Sigma_1^2}= O(n^{-1})
\hbox{ }\hbox{ and }\hbox{ }
{\Sigma_3\over \Sigma_1^2}= O(n^{-1}) \leqno(5.7)
$$
where the constants in $O(\cdot)$ depend only on $c_1,c_2$ and $\alpha$. 

\smallskip

It follows from formula (4.3) and the proof of Corollary 3.3 that for any fixed ${1\over 2}< \delta <1$
$$
\max_{\lambda_1\le \delta n}r^\lambda_2=
\max_{\lambda_1\le \delta n} {1\over {n\choose 2}}
\left[\sum\limits_i {\lambda_i\choose 2}-\sum\limits_j{\lambda_j'\choose 2}\right]\le 
$$
$$
\le {1\over {n\choose 2}}\cdot
\left[{\delta n\choose 2}+{(1-\delta)n\choose 2}\right]
\le \delta^2+(1-\delta)^2,
$$
where the maximum is taken over all partitions $\lambda$ of $n$ with $\lambda_1\le \delta n$.

Hence
$$
{1\over (1-r^\lambda_2)^2}\le {1\over 2\delta(1-\delta)}. \leqno(5.8)
$$

The absolute value of a normalized character of a
finite group is bounded above by 1.
Combining this elementary fact with (5.8)
implies that there exist constants $c_1(\delta),c_2(\delta)$
independent of $n$, so that
$$
c_1(\delta)\le {1-r^\lambda_{22}\over (1-r^\lambda_2)^2}\le 
c_2(\delta)\hbox{ }\hbox{ and }\hbox{ }
 c_1(\delta)\le {4r^\lambda_3-3r^\lambda_{22}-1\over 3(1-r^\lambda_2)^2} \le c_2(\delta). \leqno(5.9)
$$

\smallskip

To complete the proof it still remains to estimate the asymptotics 
of $r^\lambda_{22}-(r^\lambda_2)^2.$

\smallskip

Using the classical Frobenius character formula, it may be shown that
$$
r_{22}^\lambda={4\over (n-2)(n-3)}+
{1\over {n\choose 2,2,n-4}}\cdot
\left[
\left( \sum_{(i,j)\in \lambda} (j-i)\right)^2-
3\sum_{(i,j)\in \lambda} (j-i)^2 
\right] .
$$
See [I] and [Su].

Using this formula, (4.3) and (3.8) one obtains 
$$
|r^\lambda_{22}-(r^\lambda_2)^2|\le
{4\over (n-2)(n-3)}+\left[{{n\choose 2}^2\over {n\choose 2,2,n-4}}-1\right]
(r_2^\lambda)^2+   \leqno (5.10)
$$
$$
+{3\over {n\choose 2,2,n-4}}
 \sum_{(i,j)\in \lambda}(j-i)^2\le
O(n^{-2})+O(n^{-1})(r_2^\lambda)^2+O(n^{-1})=O(n^{-1}). 
$$
Substituting (5.7)-(5.10) into the right hand side of (5.6) shows that
$$
{\var_\lambda[d_f]\over E_\lambda[d_f]^2}=O(n^{-1}).
$$
This completes the proof.

\qed

\medskip 

\section{\bf Exponents in Hecke Algebras}

 Surprisingly, the expectations which appear in Proposition 3.1 and Theorem 4.1 turn out to be the exponents of $q$ in eigenvalues of irreducible Hecke algebra representations.
In particular, $E_\lambda[\descent]$ and $E_\lambda[\major]$ are the exponents of $q$ for the Hecke
 algebra elements corresponding to a Coxeter element and the longest
element of $S_n$, respectively. This follows from a well-known result of Benson and Curtis [BC].

\medskip

The Hecke algebra ${\cal H}_n(q)$ of type $A$ is the algebra over 
$\bbf := \bbc(q^{1\over 2})$
generated by $n-1$ generators  $T_1,\dots,T_{n-1}$, 
satisfying the Moore-Coxeter relations 
$$
T_iT_{i+1}T_i=T_{i+1}T_iT_{i+1} 
 \hbox{ }\hbox{ }\hbox{ }\hbox{ } (1\le i< n-1)
$$
and
$$
T_iT_j=T_jT_i  \hbox{ }\hbox{ }\hbox{ }\hbox{ if }  |i-j|>1.
$$
 as well as 
the following ``deformed involution'' relation:
$$
T_i^2= (1-q)T_i+q \qquad (1\le i<n).
$$
Note that the third relation is slightly non-standard. This is done
in order to have a more elegant formulation of Proposition 6.1. 

\medskip

Let $w_0:=(1,n)(2,n-1)\cdots$ be the longest permutation
in the symmetric group $S_n$,
and $c_n:=(12\dots n)$ be a full cycle (also known as 
a Coxeter element).

\medskip

\begin{note}
\noindent{\bf Proposition 6.1}
\begin{itemize}

\item[(i)] The eigenvalues of $T_{w_0}$ in the irreducible representation
indexed by $\lambda$ are $\pm q^{E_\lambda[\major]}$.

\item[(ii)] The eigenvalues of $T_{c_n}$ in the irreducible representation
indexed by $\lambda$ are 
$$
\omega^{\major(T)}q^{E_\lambda[\descent]}
$$
where $T$ ranges over all standard tableaux of shape $\lambda$,
and $\omega=e^{2\pi i\over n}$.

\end{itemize}
\end{note}

\smallskip

\noindent{\bf Proof.} For self containment of the paper we recall
the proof of
 [BC, corrections and additions]. See also [GM, Lemma 4.1].
It is well known that $T_{w_0}^2$ lies in the center of ${\cal H}_n(q)$.
Hence, denoting by $\rho^\lambda$
the irreducible representation of ${\cal H}_n(q)$ indexed by $\lambda$,
$\rho^\lambda(T_{w_0}^2)$ is a
scalar operator.

On the other hand for each generator $T_i, 1\le i < n$, $\rho^\lambda(T_i)$ has
two eigenvalues: $1$ and $-q$, with multiplicities
${1\over 2}(f^\lambda+\chi^\lambda_2)$ and 
${1\over 2}(f^\lambda-\chi^\lambda_2)$ respectively.
Here $f^\lambda$ is the degree of $\rho^\lambda$
and $\chi^\lambda_2=\hbox{tr}\rho^\lambda(T_i)|_{q=1}$.
 Hence
$$
\dt \rho^\lambda(T_i)=(-q)^{{1\over 2}(f^\lambda-\chi^\lambda_2)}=
\pm q^{{1\over 2}(f^\lambda-\chi^\lambda_2)}.
$$
$T_{w_0}$ is a product of ${n\choose 2}$ generators $T_i$. Hence
$$
\dt \rho^\lambda(T_{w_0}^2)=\dt \rho^\lambda(T_i)^{n(n-1)}=q^{{n\choose 2}(f^\lambda-\chi^\lambda_2)}.
$$
This shows that the eigenvalues of the scalar operator 
$\rho^\lambda(T_{w_0}^2)$ are all equal to $q^{{n \choose 2}
{f^\lambda-\chi^\lambda_2\over f^\lambda}}$, and those of $T_{w_0}$ are
$$
\pm q^{{n\choose 2}{f^\lambda-\chi^\lambda_2\over 2f^\lambda}}=
\pm q^{{n\choose 2}{1-r^\lambda_2\over 2}}=\pm q^{E_\lambda[\major]}.
$$
The last equality follows from Proposition 3.1 and formula (4.3).

\smallskip

To prove the second part of Proposition 6.1 recall that
$$
T_{c_n}^n=T_{w_0}^2.
$$
So, the eigenvalues of $T_{c_n}$ are the complex $n$-th roots of the eigenvalues of $\rho^\lambda(T_{w_0}^2)$. Combining this fact with Theorem 4.1 gives the 
corresponding exponent of $q$.
For a calculation of the exponents of $\omega$ see [Ste].

\qed

\smallskip

Proposition 6.1 may be useful in the study of probabilistic interpretations of the Hecke algebra, and of related asymmetric random walks.


\begin{thebibliography}{CF3}

\bibitem[AR]{AR} R. M. Adin and Y. Roichman,
{\it On random Young tableaux (extended abstract)},
In: Paul Erd\H os and his Mathematics - Research Communications,
J\'anos Bolyai Mathematical Society, Budapest, 1999, pp. 4--6.


\bibitem[AS]{AS} N. Alon and J. H. Spencer,  {\it The Probabilistic Method},
(with an appendix by Paul Erd\H os). Wiley, New-York, 1992.

\bibitem[BC]{BC} C. T. Benson and C. W. Curtis,
{\it On the degrees and rationality of certain characters of finite Chevalley
groups}, Trans. Amer. Math. Soc. 165 (1972), 251--273;
{\it Corrections and additions}, ibid. 202 (1975), 405--406.

\bibitem[DMP]{dmp} P. Diaconis, M. McGreath and J. Pitman,
{\it Riffle shuffles, cycles, and descents.}
Combinatorica 15 (1995), 11--29.

\bibitem[Fe]{Fe} W. Feller, {\it An Introduction to Probability Theory
and its Applications, Vol. I}, Wiley, New-York, 1970.

\bibitem[F]{F} 
D.\ Foata, 
{\it On the Netto inversion number of a sequence.}
Proc.\ Amer.\ Math.\ Soc.~19 (1968), 236--240.

\bibitem[FS]{FS} 
D.\ Foata and M.\ P.\ Sch\"utzenberger,
{\it Major index and inversion number of permutations.} 
Math.\ Nachr.~83 (1978), 143--159.

\bibitem[Fu]{Fu} J. Fulman, {\it The distribution of descents in fixed
conjugacy classes of the symmetric groups}, J. Combin. Theory Ser. A 84
(1998), 171--180.

\bibitem[GG]{GG}
A.\ Garsia and I.\ Gessel,
{\it Permutation statistics and partitions.}
Adv.\ in Math.~31 (1979), 288--305. 

\bibitem[GM]{GM} M. Geck and J. Michel, {\it `Good' elements of finite Coxeter groups and representations of Iwahori-Hecke algebras},
Proc. London Math. Soc. 74 (1997), 275--305.

\bibitem[H]{h} P. A. H\"ast\"o, {\it On descents in standard Young tableaux},
preprint, July 2000. 

\bibitem[I]{I} R. E. Ingram, {\it Some characters of the symmetric group}, Proc. Amer. Math. Soc. 1 (1950), 358--369.

\bibitem[K]{k} D. Knuth, {\it The Art of Computer Programming,
Volume 3 : Sorting and Searching}, Adisson-Wesley,
Reading MA, 1973.

\bibitem[M]{M} 
P.\ A.\ MacMahon, 
{\it Combinatory Analysis, Volumes I-II}. 
Cambridge Univ.\ Press, London/New-York, 1916. 
(Reprinted by Chelsea, New-York, 1960.)

\bibitem[Re]{re} C.\ Reutenauer, 
Free Lie Algebras. 
London Mathematical Soc.\ Monographs, New Series 7, Oxford Univ.\ Press, 1993.

\bibitem[Ro1]{Ro} Y. Roichman, {\it A recursive rule for Kazhdan-Lusztig characters.} Adv. in Math. 129 (1997), 25--45.

\bibitem[Ro2]{ro} Y. Roichman, {\it On Permutation Statistics and  Hecke Algebra Characters},
In: Combinatorial Methods in Representation Theory,
Adv. Pure Math., Math. Soc. Japan, to appear.

\bibitem[Sa]{Sa} B. E. Sagan, {\it The Symmetric Group. Representations,
Combinatorial Algorithms, and Symmetric Functions.} Wadsworth \& Brooks/Cole,
CA, 1991.

\bibitem[St]{St} R. P. Stanley, {\it Enumerative Combinatorics, Volume II.}
Cambridge Univ. Press, Cambridge, 1999. 

\bibitem[Ste]{Ste} J. Stembridge, {\it On the eigenvalues of representations of reflection
groups and wreath products.} Pacific J. Math. 140 (1989), 359--396.

\bibitem[Su]{Su}  M. Suzuki, {\it The values of irreducible characters of the symmetric group}, The Arcata Conference on Representations of Finite Groups,
Amer. Math. Soc. Proceedings of Symposia in Pure Mathematics, Vol. 47 - Part 2 (1987), 317--319.

\end{thebibliography}
\end{document}